\documentclass{article} \usepackage{graphicx,amsmath,amssymb,amsthm}

\theoremstyle{plane}
\newtheorem*{thm}{Theorem}
\newtheorem*{lem}{Lemma}

\DeclareMathOperator{\Isom}{Isom}
\DeclareMathOperator{\Aut}{Aut}

\newenvironment{abenumerate}{%
	 \begin{enumerate}%
	}{\end{enumerate}%
}

\title{A classification of hyperbolic manifolds related to the three-dimensional torus}
\author{Akira {\sc Ushijima} \vspace{\medskipamount}\\
Mathematics Institute\\
University of Warwick\\
Coventry CV4 7AL, United Kingdom\\
{\tt ushijima@maths.warwick.ac.uk}}
\date{\empty}
\begin{document} \maketitle
\begin{abstract}
L.~Paoluzzi constructed a family of compact 
orientable three-dimensional 
hyperbolic manifolds with totally geodesic boundary,
which were, by construction, closely related to the three-dimensional torus.
This paper gives their complete classification up to isometry,
and also their isometry groups.
The key tool is the so-called canonical decomposition of hyperbolic manifolds.
\vspace{\medskipamount}
\begin{flushleft}
		\vspace{\smallskipamount}
\end{flushleft}
\end{abstract}
%
%
%
%
%
%
%
%
\section{Introduction}
In \cite{pa} Paoluzzi constructed 
a family $\left\{ M_{n,k} \right\}$ of
compact orientable three-dimensional manifolds
with boundary, indexed by integers
$n$ and $k$ with $n \geq 1$ and $0 \leq k \leq n-1$.
Though the construction is closely related to the three-dimensional torus ${\bf T}^3$,
one of the typical non-hyperbolic manifolds
(and $M_{3,1}$ is indeed, by definition, ${\bf T}^3$ minus a small open ball),
most of them are shown to admit a hyperbolic structure 
with totally geodesic boundary.
The key property for the construction
is the symmetry of ${\bf T}^3$ of order three induced by the diagonal of the cube 
as its fundamental polyhedron.
It should be noted that a similar construction appears to give hyperbolic spatial graphs 
from torus knots in the three sphere (see \cite{us_HypGraph}).

Another interesting property of $M_{n,k}$ is that 
its fundamental group has cyclic presentation.
A group is said to have {\em cyclic presentation}\/
if it is isomorphic to
$\left< x_0 , x_1 , \ldots , x_{n-1} 
	\left| \, w = \theta (w) = \cdots = \theta^{n-1} (w) \right. \right>$
for some $n \in \mathbb{N}, w$ and $\theta$,
where $w$ is a word in the free group
$F_n := \left< \left. x_0 , x_1 , \ldots , x_{n-1} \right| \, \right>$
and $\theta$ is the automorphism of $F_n$
defined by $\theta (x_i) := x_{i+1}$ for all $i = 0 , 1 , \ldots , n-1$,
where all indices are considered modulo $n$.
It has already been shown in \cite{st} that 
it is impossible to find an algorithmic solution to
whether or not a finite presentation of a group is
isomorphic to the fundamental group of a three-dimensional manifold.
Furthermore, cyclically presented groups contain well-known and well studied groups
such as Fibonacci groups and Sieradski groups,
which it is known can be regarded as
fundamental groups of some three-dimensional hyperbolic manifolds
(see, for example, \cite{bv,cm} and references therein
for 
related topics).
So it is worthwhile (cf.\ \cite{chr,du1}) to investigate 
three-dimensional hyperbolic manifolds 
whose fundamental groups are isomorphic to some other cyclically presented groups,
and $M_{n,k}$ is indeed such a manifold.

\vspace{\bigskipamount}

\begin{figure}[ht]
        \begin{center}
        \includegraphics[width=50mm,clip]{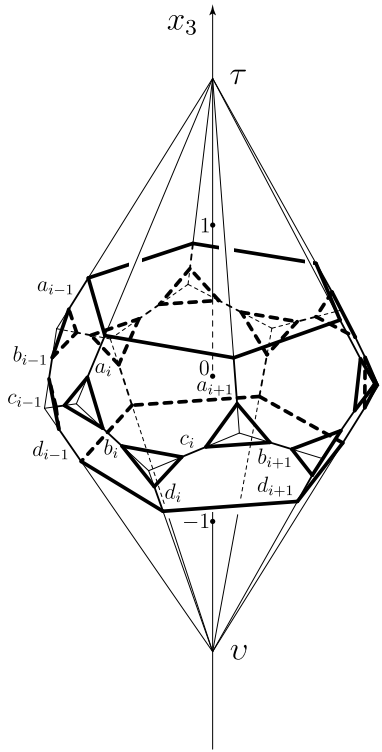}
        \end{center}
        \caption{The polyhedra $\mathcal{CP}_6$ and $\mathcal{P}_6$}
        \label{fig_cpol}
\end{figure}

We here recall a combinatorial definition of the manifold $M_{n,k}$ given in \cite{pa}.
For a positive integer $n$, 
let $\mathcal{CP}_n$ be a dual polyhedron of a
regular $n$-gonal anti-prism 
(see Figure~\ref{fig_cpol}); 
it is obtained by taking the quotient of a cube by
the rotation of order three about its diagonal,
and then taking the branched covering of the 
``$1/3$-cube" along the image of the diagonal of order $n$. 
So its boundary consists of $2 \, n$ quadrilateral faces.
We could say that $\mathcal{CP}_n$ is a natural generalization of a cube.
Indeed, it plays a key role in {\em generalized cubing}\/ 
of three-dimensional manifolds (see \cite[p.\ 132]{ar}).
We then denote by $\mathcal{P}_n$ the polyhedron obtained
by truncating the neighborhoods of all vertices of $\mathcal{CP}_n$ 
in such a small way that the faces obtained by the truncation do not intersect each other
(see Figure~\ref{fig_cpol} again). 
We call a face of $\mathcal{P}_n$ obtained by a truncation 
{\em external\/} ({\em internal\/} otherwise). 
Similarly we call an edge surrounding an external
face {\em external\/} ({\em internal\/} otherwise). 
%
%
%
%
%
%

\begin{figure}[ht]
        \begin{center}
        \includegraphics[width=100mm,clip]{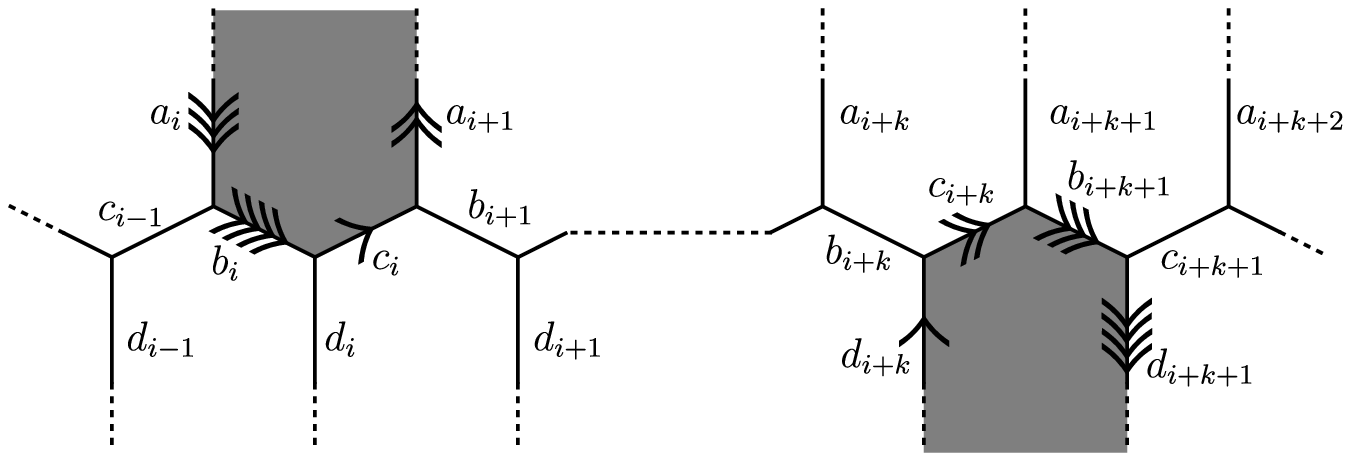}
        \end{center}
        \caption{The identification rule of step $k$}
        \label{fig_idrule}
\end{figure}

We next give an identification rule among the
internal faces of $\mathcal{P}_n$. 
Figure~\ref{fig_idrule} shows the boundary of 
$\mathcal{CP}_n$ flattened out on the plane,
where the top (resp.\ bottom) cone point is 
at positive (resp.\ negative) infinity. 
So we can imagine that the solid occupies the
half-space behind the page. 
We fix an integer $k$ with $0 \leq k \leq n-1$. 
Then, for each integer $i$ with $0 \leq i \leq n-1$, we identify the face 
$c_i \, a_{i+1} \, a_i \, b_i$ with the face 
$d_{i+k} \, c_{i+k} \, b_{i+k+1} \, d_{i+k+1}$ 
in this order, where indices are always considered modulo $n$. 
We call this identification rule the 
{\em identification rule of step $k$}. 
We here note that, when $n=3$, $\mathcal{CP}_3$ is
simply a cube, and then 
the identification rule of step $1$ 
gives the three-dimensional torus ${\bf T}^3$.
This identification rule naturally
induces the one among the internal
faces of $\mathcal{P}_n$, 
and we denote a topological manifold of this quotient space by $M_{n,k}$. 
%
%
%
%
%
%

It should be noted 
that there are two different behaviors of $M_{n,k}$ depending on $n$. 
When $n \not \equiv 0$ modulo $3$,
all internal edges of $\mathcal{P}_n$ are
glued together, and $\partial M_{n,k}$ is
an orientable surface of genus $n - 1$. 
On the other hand, when $n \equiv 0$ modulo $3$,
the internal edges are identified in three groups of $4 \, n / 3$ each, 
and $\partial M_{n,k}$ is an orientable surface of genus $n - 3$. 

\vspace{\bigskipamount}

\begin{figure}[ht]
        \begin{center}
        \includegraphics[clip]{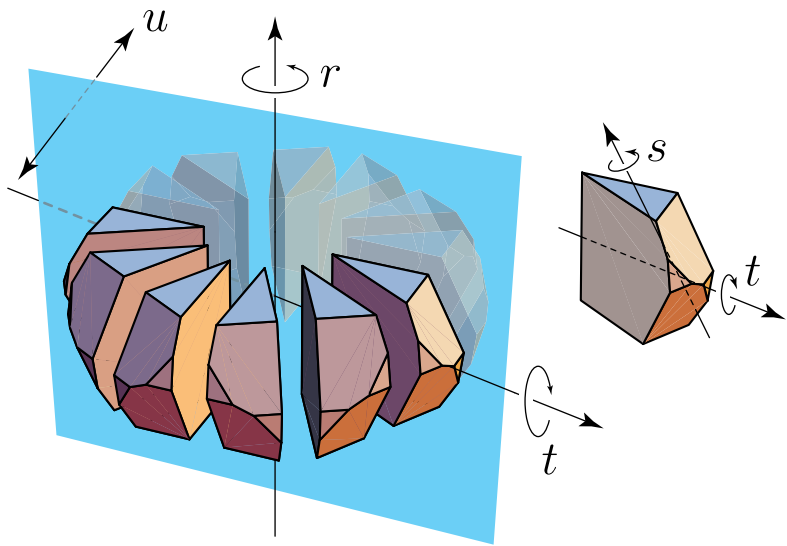}
        \end{center}
        \caption{Generators of $\Isom M_{n,k}$}
        \label{fig_sym}
\end{figure}

It is shown in \cite{pa}
that each $M_{n,k}$ is a cyclic branched covering of 
$\left( {\bf S}^2 \times {\bf S}^1 \right) - {\bf B}^3$ along two arcs,
and that it admits a hyperbolic structure with totally geodesic boundary 
when $n \geq 4$, except $n = 6$, and any $k$ with $0 \leq k \leq n-1$.
She also determined their fundamental groups and the first homology groups,
and using them she obtained
a sufficient condition for two manifolds to be isometric.
Indeed she proved that $M_{n,k}$ and $M_{n',k'}$ are isometric
if $n' = n$ and $k' = k$ or $n - k - 1$.
In this paper, we complete her results.
Namely we prove the following theorem,
which improves her results to say that 
$M_{n,k}$ admits a hyperbolic structure even when $n=6$, 
and that her sufficient condition is enough to classify the manifolds.
We also obtain their isometry groups as a by-product of the classification.
The key tool for the proof is the canonical decomposition of $M_{n,k}$,
which will be explained in the proof given in the next section.
%
%
%
%
%
%
%
%
\begin{thm}
\begin{abenumerate}
\item \label{it_str}
The manifolds $M_{n,k}$ are hyperbolic with totally geodesic boundary
when $(n,k)$ satisfies $n \geq 4$ and $0 \leq k \leq n-1$.
\item \label{it_class}
Two manifolds\/ $M_{n,k}$ and\/ $M_{n',k'}$ are isometric if and only if\/
$n' = n$ and\/ $k' = k$ or\/ $n - k - 1$. \label{thm_isom_homeo}
\item \label{it_isom}
The isometry group\/
$\Isom M_{n,k}$
of $M_{n,k}$ is as follows:
\begin{itemize}
\item Suppose $n \not\equiv 0$ modulo $3$. Then
\begin{itemize}
\item $\Isom M_{n,k} \cong \left< t, u \left| \, t^2 = u^2 = (ut)^{2 \, n} = 1 \right. \right>$
when $n$ is odd and $k = (n-1)/2$,
\item $\Isom M_{n,k} \cong \left< r, t \left| \, r^n = t^2 = (tr)^2 = 1 \right. \right>$
otherwise.
\end{itemize}
\item Suppose $n \equiv 0$ and 
$k \not\equiv 1$ modulo $3$. Then\/
$\Isom M_{n,k} \cong \left< r, t \left| \, 
r^n = t^2 = (tr)^2 = 1 \right. \right>$.\sloppy
\item Suppose $n \equiv 0$ and $k \equiv 1$ modulo $3$, i.e., 
$n = 3 \, m$ for some $m \geq 2$ and 
$k = 3 \, l + 1$ for some $l \geq 0$. Then 
\begin{itemize}
\item $\Isom M_{n,k} \cong 
\left< s, t, u \left| \, \begin{array}{l}s^2 = t^2 = u^2 = (st)^2 = (ut)^6 = 1, \\
susus = tutut \end{array} \right. \right>$
when $m$ is odd and $l = (m-1)/2$,
\item $\Isom M_{n,k} \cong \left< r, s, t \left| \, 
\begin{array}{l}
r^{3 \, m} = s^2 = t^2 = (tr)^2 = (st)^2 = 1,\\
sr^3 = r^3s, (str)^3 = r^{3 \left( m - 2 \, l - 2 \right)} 
\end{array} \right. \right>$ otherwise.
\end{itemize}
\end{itemize}
\end{abenumerate}
In these presentations, $r$ is induced by the\/
$2 \, \pi /n$-rotation along $x_3$ axis, 
$t$ is induced by the $\pi$-rotation along the axis joining the origin and 
the midpoint of the edge obtained from $b_0$, 
$u$ is the reflection of the plane 
containing $x_3$ axis and $d_0$, 
$s$ is induced by the $\pi$-rotation along the axis 
joining the midpoint of the edge obtained from $a_0$ 
and the midpoint of the edge obtained from $d_0$ (see Figure~\ref{fig_sym}).
\end{thm}
%
%
%
%
%
%
%
%
%
%
%
%
\subsubsection*{Acknowledgement}
The author would like to thank 
Professor Makoto Sakuma and Dan Goodman
for his helpful comments and advice.
The author would also like to express 
his sincere gratitude to Professor Susumu Hirose
for useful discussions on the proof of part \eqref{it_isom} of the theorem,
especially on the combinatorial group theory.
%
%
%
%
%
%
%
%
\section{Proof of the theorem}
%
%
%
%
%
%
\subsection*{Fundamental polyhedron of $M_{n,k}$}
First we show \eqref{it_str}.
Though Paoluzzi proved in \cite{pa} (most of) this result using Andreev's theorem,
here we determine the fundamental polyhedron of $M_{n,k}$ explicitly,
for use later in the proof.
The outline of this part is as same as of \cite[Theorem~2.1]{us_canon},
and we omit basic results about hyperbolic geometry.

Let $\mathbf{P}_1^{3} := \left\{ \left. \boldsymbol{x} 
\in \mathbb{E}^{1,3} \, \right| x_0 = 1 \, \right\}$
be an affine hyperplane in the four-dimensional Minkowsky space
$\mathbb{E}^{1,3} = 
(\mathbb{E}^4, \left \langle \, \cdot , \cdot \, \right \rangle)$.
We regard $\mathbf{P}_1^{3}$ as $\mathbb{E}^3$,
and we also regard the open unit ball, say $\mathbf{B}^3$, centred at the origin 
as the projective ball model of three-dimensional hyperbolic space, 
in the usual way. That is, the radial projection from
$\left\{ \left. 
\boldsymbol{x} \in \mathbb{E}^{1,3} \, \right|
\left \langle \, \boldsymbol{x}, 
\boldsymbol{x} \, \right \rangle 
= - \, 1 \text{ and } x_0 > 0 \, \right\}$,
the so-called hyperboloid model of the hyperbolic space.
Then $H_S := \left\{ \left. \mbox{\boldmath $x$} 
\in {\Bbb E}^{1,3} \, \right|
\left \langle \, \mbox{\boldmath $x$}, 
\mbox{\boldmath $x$} \, \right \rangle  
= 1 \, \right\}$, which corresponds to the set of half-spaces in the hyperbolic space,
projects to the exterior of the closed unit sphere $\overline{\mathbf{B}^3}$.

\begin{figure}[ht]
        \begin{center}
        \includegraphics[width=110mm,clip]{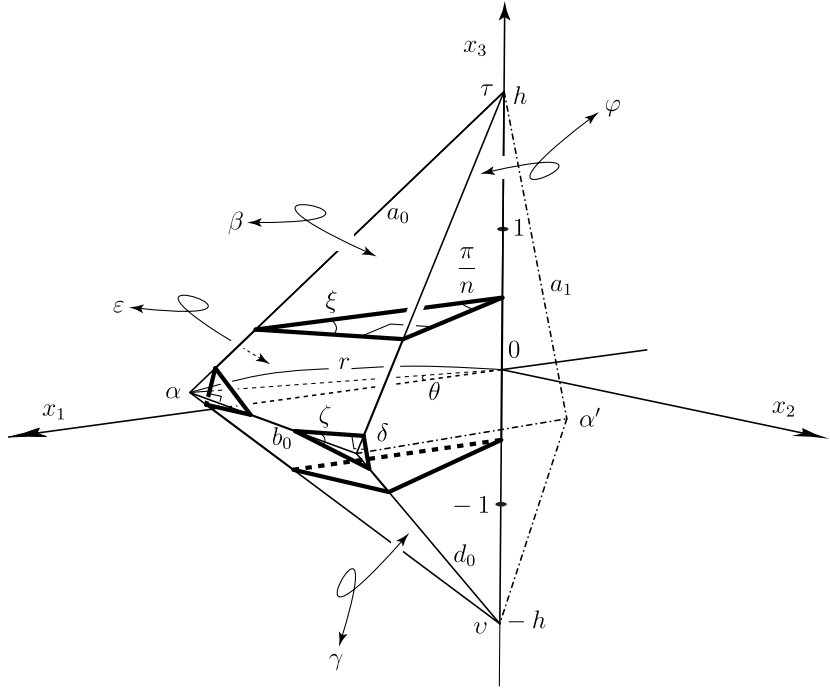}
        \end{center}
        \caption{A tetrahedron ${\cal Q}_n$}
        \label{fig_piece}
\end{figure}

Suppose ${\cal CP}_n$ is a Euclidean symmetric polyhedron
which is put in the centre of $\mathbf{P}_1^{3}$ as in Figure~\ref{fig_cpol}.
Since we are to construct the fundamental polyhedra of the hyperbolic manifolds,
we can also assume that ${\cal CP}_n$ is symmetric and its centre
coincides with the origin of $\mathbf{P}_1^{3}$,
and we focus on a fundamental tetrahedron, say ${\cal Q}_n$,
with respect to this symmetry (see Figure~\ref{fig_piece}).
%
%
%
%
%
%

Let $\tau$ (resp.\ $\upsilon$) be the top (resp.\ bottom) cone point,
and let $\alpha$ (resp.\ $\delta$) be the terminal point of the edge 
$a_0$ (resp.\ $d_0$) opposite to $\tau$ (resp.\ $\upsilon$).
These vertices are poles of truncation planes,
so they must lie outside of $\overline{\mathbf{B}^3}$.
We denote by $\widehat{\cdot}$ their lifts to $H_S$.
Now the coordinates of 
$\widehat{\tau}$, $\widehat{\upsilon}$, $\widehat{\alpha}$ and $\widehat{\delta}$
can be expressed by $h$, $r$ and $\theta$ as follows:
\begin{equation*}
\widehat{\tau} = \frac{1}{\sqrt{h^2 - 1}} \left( 1, 0, 0, h \right) \, , \,
\widehat{\upsilon} = \widehat{\tau} \, N \, , \,
\widehat{\alpha} = 
	\frac{1}{\sqrt{r^2 - 1}} \left( 1, r \cos \theta, 0, r \sin \theta \right) \, , \,
\widehat{\delta} = \widehat{\alpha} \, N  \, .
\end{equation*}
Here $h$ (resp.\ $r$) is the Euclidean distance 
between $\tau$ (resp.\ $\alpha$) and the origin,
$\theta$ is the Euclidean angle between the $x_1$-axis
and the segment joining $\alpha$ and the origin,
and $N$ is a $4 \times 4$ matrix defined as follows:
\begin{equation*}
N := \begin{pmatrix}
1 & 0 & 0 & 0 \\
0 & c_n & \sqrt{1 - c_n^2} & 0 \\
0 & - \, \sqrt{1 - c_n^2} & c_n & 0 \\
0 & 0 & 0 & - \, 1 
\end{pmatrix} \text{ , where } c_n := \cos \frac{\pi}{n} \, .
\end{equation*}
Let $\widehat{\gamma}$ be the outward unit normal vector to the plane spanned by 
$\widehat{\alpha}$, $\widehat{\delta}$ and $\widehat{\upsilon}$.
Then its coordinate is:
\begin{multline*}
\widehat{\gamma} = \frac{1}{\sqrt{H_1}}
( h \, r \, \sqrt{1 - c_n^2} \, \cos \theta , 
\left( h + r \sin \theta \right) \sqrt{1 - c_n^2} , \\
\qquad \qquad \left( 1 - c_n \right) h - r \left( 1 + c_n \right) \sin \theta , 
- \, r \, \sqrt{1 - c_n^2} \, \cos \theta ) \, ,
\end{multline*}
where $H_1 := 
r^2 \left( 1 + c_n \right) \left\{ 1 + c_n + \left( 1 - c_n \right) h^2 \right\} \sin^2 \theta
+ r^2 \left( 1 - c_n^2 \right) \left( 1 - h^2 \right) + 2 \left( 1 - c_n \right) h^2$.
Then $\widehat{\beta} := \widehat{\gamma} \, N$ is the
outward unit normal vector to the plane spanned by
$\widehat{\delta} = \widehat{\alpha} \, N$,
$\widehat{\alpha'} = \widehat{\delta} \, N$ and
$\widehat{\tau} = \widehat{\upsilon} \, N$.

By the definition of $\mathcal{CP}_n$,
the four vertices $\tau$, $\alpha$, $\delta$ and $\alpha'$ must be in the same plane.
This condition is equivalent to $\langle \, \widehat{\alpha}, \widehat{\beta} \, \rangle = 0$.
Thus we obtain the following relation:
\begin{equation*}
\sin \theta = \frac{\left( 1 - c_n \right) h}{\left( 1 + c_n \right) r} \, .
\end{equation*}

For isometric gluing, the hyperbolic lengths of the edges
obtained from $a_0$ and $b_0$ by truncation must be equal.
This condition is equivalent to 
$\langle \, \widehat{\tau}, 
\widehat{\alpha} \, \rangle 
= \langle \, \widehat{\alpha}, 
\widehat{\delta} \, \rangle$,
and we thus obtain the following relation:
\begin{equation*}
r = \frac{\sqrt{2 \, c_n - 1 + (1 - c_n)^2 \, h^2}} {c_n} \, .
\end{equation*}

Let $\widehat{\varepsilon}$ (resp.\ $\widehat{\varphi}$) be
 the outward unit normal vector to the plane spanned by 
$\widehat{\tau}$, $\widehat{\alpha}$ and 
$\widehat{\upsilon}$ (resp.\ $\widehat{\tau}$, 
$\widehat{\delta}$ and $\widehat{\upsilon}$).
Then $\widehat{\varepsilon} = \left( 0, 0, - \, 1, 0 \right)$ and
$\widehat{\varphi} = \widehat{\varepsilon}$.
Now the dihedral angles of $\mathcal{Q}_n$ can be calculated by using
$\widehat{\beta}$, $\widehat{\gamma}$ and $\widehat{\varepsilon}$.
If we denote by $\xi$ (resp.\ $\zeta$)
the dihedral angle between polar planes of $\widehat{\beta}$ and $\widehat{\varepsilon}$
(resp.\ $\widehat{\beta}$ and $\widehat{\gamma}$),
then $\xi = \arccos ( - \, \langle \, \widehat{\beta}, \widehat{\varepsilon} \, \rangle )$ 
and $\zeta = \arccos ( - \, \langle \, \widehat{\beta}, \widehat{\gamma} \, \rangle)$.
Since the dihedral angles of edges glued together must sum up to $2 \, \pi$,
we have the third condition on the parameters $h$, $r$ and $\theta$
via a relation with respect to $\xi$ and $\zeta$.
Since there are two kinds of gluing patterns,
as mentioned in the middle of the introduction,
the condition separates into two cases.

We first consider the case where 
$n \not \equiv 0$ modulo $3$. 
In this case all internal edges are identified together with one edge. 
So we obtain the condition 
$2 \, n \left( 2 \, \xi + \zeta \right) = 2 \, \pi$. 
By straightforward calculation, it is equivalent to the following one:
\begin{equation*}
h = \frac{\sqrt{1 + c_n} \, \sqrt{\sqrt{8 \, c_n^2 + 1} - 3 \, c_n}}{1 - c_n} \, .
\end{equation*}
Now we can easily check that, 
when $n \geq 4$ with $n \not \equiv 0$ modulo $3$, 
$h, r \geq 1$ and the truncation planes are ultraparallel.
Thus we have obtained a suitable realization of $\mathcal{P}_n$
as a hyperbolic polyhedron in this case.

We next consider the case where $n \equiv 0$ modulo $3$. 
In this case the internal edges are identified in three groups of $4 \, n / 3$ each. 
So we obtain the condition 
$\left. 2 \, n \left( 2 \, \xi + \zeta \right) \right/ 3 = 2 \, \pi$. 
In the same fashion as above, we have
\begin{equation*}
h = \frac{\sqrt{1 + c_n} \, \sqrt{2 \, c_n - 1}}{\sqrt{1 - c_n} \, \sqrt{2 \, c_n + 1}} \, ,
\end{equation*}
and we can see that this also gives a suitable realization of $\mathcal{P}_n$.

Finally by Poincar\'e's theorem on fundamental polyhedra (see \cite{ma}), 
we can say that the hyperbolic polyhedron ${\cal P}_n$ is
a fundamental polyhedron of $M_{n,k}$ when $n \geq 4$.
Thus \eqref{it_str} has been proved.
It should be emphasized that $M_{6,k}$ admits a hyperbolic structure, 
though Paoluzzi did not mention it in \cite{pa}. 
%
%
%
%
%
%
\subsection*{Canonical decomposition of $M_{n,k}$}
For the proof of \eqref{it_class} and \eqref{it_isom},
we use the {\em canonical decomposition}\/ of $M_{n,k}$.
For details about the canonical decomposition of 
compact hyperbolic manifolds with totally geodesic boundary, see \cite{ko1,ko2}.
Two remarkable applications about the canonical decomposition are: 
(i) two hyperbolic manifolds are homeomorphic 
(or equivalently, by Mostow's rigidity theorem, 
isometric) if and only if their canonical
decompositions are equivalent, that is, 
they are decomposed into same pieces, 
and the identification rule is also the same;
(ii) all combinatorial automorphisms of the canonical decomposition 
form the isometry group of the manifold.
So, to prove \eqref{it_class} and \eqref{it_isom},
we need to determine the canonical decomposition of $M_{n,k}$,
and the result is as follows:
%
%
%
%
%
%
%
%
\begin{lem}
The canonical decomposition\/ $\mathcal{D}_{n,k}$ of\/ $M_{n,k}$, 
where\/ $n \geq 4$ and\/ $0 \leq k \leq n-1$,  
consists of\/ $2 \, n$ truncated tetrahedra obtained by cutting\/ $\mathcal{P}_n$ 
into pieces like segments of an orange (as in Figure~\ref{fig_sym}). 
\end{lem}
\begin{proof}
If the candidate decomposition consists of truncated tetrahedra then,
by the definition of the canonical decomposition,
checking that all the tilts are negative is enough to show that the decomposition is canonical
(see \cite{us_tilt} for the definition of the tilts of truncated polyhedra).
By symmetry of $\mathcal{P}_n$, all we have to show is that all the tilts of 
the truncated tetrahedron $\mathcal{Q}_n$ are negative (see Figure~\ref{fig_piece} again). 

Let $t_{\beta}$, $t_{\gamma}$, $t_{\varepsilon}$, 
$t_{\varphi}$ be the tilts of internal faces of  
$\mathcal{Q}_n$ whose normal vectors are
$\widehat{\beta}$, 
$\widehat{\gamma}$, 
$\widehat{\varepsilon}$, 
$\widehat{\varphi}$ respectively. 
Then, by \cite[Theorem~3.3]{us_tilt}, the following equation holds:
\begin{equation*}
\begin{pmatrix}
t_{\beta} \\ 
t_{\gamma} \\ 
t_{\varepsilon} \\ 
t_{\varphi}
\end{pmatrix}
= \begin{pmatrix}
1 & \langle \, \widehat{\beta}, 
\widehat{\gamma} \, \rangle  & 
\langle \, \widehat{\beta}, 
\widehat{\varepsilon} \, \rangle  & 
\langle \, \widehat{\beta}, 
\widehat{\varphi} \, \rangle  \\  
\langle \, \widehat{\gamma}, 
\widehat{\beta} \, \rangle  & 1 &
\langle \, \widehat{\gamma}, 
\widehat{\varepsilon} \, \rangle  & 
\langle \, \widehat{\gamma}, 
\widehat{\varphi} \, \rangle  \\ 
\langle \, \widehat{\varepsilon}, 
\widehat{\beta} \, \rangle  & 
\langle \, \widehat{\varepsilon}, 
\widehat{\gamma} \, \rangle  & 1 &
\langle \, \widehat{\varepsilon}, 
\widehat{\varphi} \, \rangle  \\ 
\langle \, \widehat{\varphi}, 
\widehat{\beta} \, \rangle  & 
\langle \, \widehat{\varphi}, 
\widehat{\gamma} \, \rangle  & 
\langle \, \widehat{\varphi}, 
\widehat{\varepsilon} \, \rangle & 1
\end{pmatrix}
\begin{pmatrix}
- \, \langle \, \widehat{\beta}, 
\widehat{\tau} \, \rangle^{- \, 1} \\ 
- \, \langle \, \widehat{\gamma}, 
\widehat{\upsilon} \, \rangle^{- \, 1} \\ 
- \, \langle \, \widehat{\varepsilon}, 
\widehat{\alpha} \, \rangle^{- \, 1} \\ 
- \, \langle \, \widehat{\varphi}, 
\widehat{\delta} \, \rangle^{- \, 1}
\end{pmatrix} \, .
\end{equation*}
By straightforward calculation, we have
\begin{align*}
t_{\beta}& = t_{\gamma} =
- \, h \sqrt{\frac{H_3}{H_2}} \left\{ 
( 1 - c_n )^2 ( 2 c_n + 1 ) h^2 + \sqrt{1 - c_n} ( 1 + c_n ) ( 2 c_n - \sqrt{1 - c_n} ) 
\right\} , \\
t_{\varepsilon}& = t_{\varphi} =
- \sqrt{H_3} \, \sqrt{1 - c_n} \, \sqrt{1 + c_n} \, 
\left( c_n - \sqrt{1 - c_n} \right) \, ,
\end{align*}
where
\begin{align*}
H_2& := \left( 1 + c_n \right) 
\left\{ \left( 1 - c_n \right) h^2 + 
\left( 1 + c_n \right) 
\left( 2 \, c_n - 1 \right) \right\} \\
& \qquad + h^2 \left( 1 - c_n \right) 
\left\{ 1 + c_n - \left( 1 - c_n \right) 
\left( 2 \, c_n + 1 \right) h^2 \right\} \, , \\
H_3& := \left(h^2 - 1\right)
\left\{ \left( 1 + c_n \right)^2 \left( 2 \, c_n - 1 \right) + 
\left( 1 - c_n \right)^2 \left( 2 \, c_n + 1 \right) h^2 \right\}^{-1} \, .
\end{align*}
Since $n \geq 4$, we have 
$2 \, c_n - \sqrt{1 - c_n} > c_n - \sqrt{1 - c_n} > 0$. 
Thus we obtain $t_{\beta} , t_{\gamma} , t_{\varepsilon} , t_{\varphi} < 0$ 
for any $n \geq 4$, thereby proving the lemma.
\end{proof}
%
%
%
%
%
%
%
%
\subsection*{Proof of \eqref{it_class} and \eqref{it_isom}}
\begin{figure}[ht]
	\begin{center}
	\includegraphics[clip]{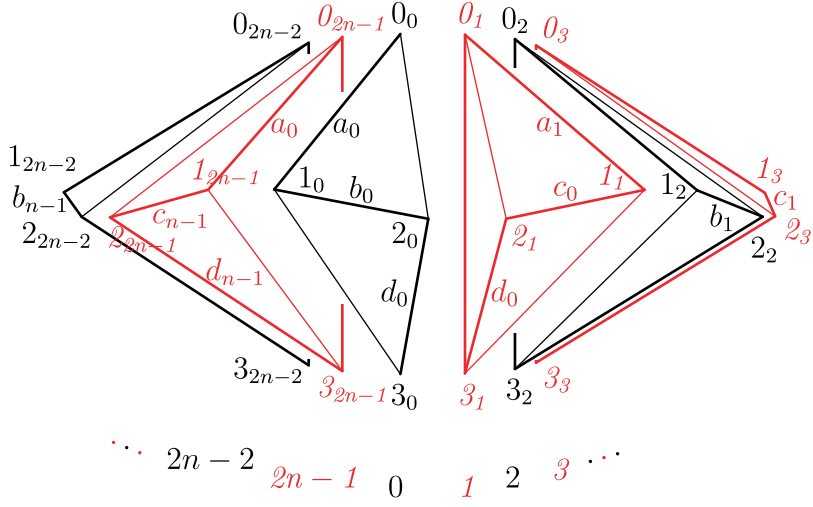}
	\end{center}
	\caption{The tetrahedral decomposition 
		$\mathcal{CD}_{n,k}$ of $\mathcal{CP}_n$}
	\label{fig_cdp}
\end{figure}

One of the most remarkable characteristics of the
canonical decomposition is that the decomposition
is invariant under the action of the isometry group of $M_{n,k}$. 
Let $\Aut \left( \mathcal{D}_{n,k} \right)$ be
the combinatorial automorphism group of the 
canonical decomposition $\mathcal{D}_{n,k}$ of $M_{n,k}$.
By Mostow's rigidity theorem together with the characteristic of the 
canonical decomposition mentioned above,
the isometry group $\Isom \left( M_{n,k} \right)$ of $M_{n,k}$ is isomorphic to
$\Aut \left( \mathcal{D}_{n,k} \right)$. 
Thus, for the study of $\Isom \left( M_{n,k} \right)$,
it is enough to investigate $\Aut \left( \mathcal{D}_{n,k} \right)$. 
Since this group is determined from the combinatorial data of $\mathcal{D}_{n,k}$, 
we focus on the (untruncated) tetrahedral decomposition,
say $\mathcal{CD}_{n,k}$, of $\mathcal{CP}_n$ 
induced by the canonical decomposition of $M_{n,k}$ 
with the identification rule of step $k$ defined in the introduction, 
and its combinatorial isomorphisms (see Figure~\ref{fig_cdp}). 
%
%
%
%
%
%

Let $\varphi$ be such an isomorphism.
Then there are two possibilities; 
the first is that $\varphi$ is indeed a combinatorial automorphism of $\mathcal{D}_{n,k}$. 
In this case $\varphi$ is an element of $\Aut \left( \mathcal{D}_{n,k} \right)$,
and so is an element of $\Isom \left( M_{n,k} \right)$. 
The second is that $\varphi$ is not 
a combinatorial automorphism of $\mathcal{D}_{n,k}$, 
but $\varphi \left( \mathcal{D}_{n,k} \right)$ 
becomes $\mathcal{D}_{n,k'}$ for some $k'$. 
In this case $\varphi$ induces an isometry between $M_{n,k}$ and $M_{n,k'}$. 
We will find all isomorphisms of the first case,
and these form all of $\Aut \left( \mathcal{D}_{n,k} \right)$,
which finishes the proof of \eqref{it_isom}.
We will find all isomorphisms of the second case,
and this gives us all isometries $M_{n,k} \cong M_{n,k'}$,
which finishes the proof of \eqref{it_class}.
So our task is to find all $\varphi$'s of the first and second cases. 

We name pieces and edges of $\mathcal{CD}_{n,k}$ as follows: 
a piece including $b_i$ (resp.\ $c_i$) is named $2 i$ (resp.\ $2 i + 1$). 
In a piece $j$, the vertex obtained from $\tau$ 
(resp.\ $\upsilon$) is named $0_j$ (resp.\ $3_j$), 
and the terminal point of $a_*$ (resp.\ $d_*$) 
other than $0_j$ (resp.\ $3_j$) is named $1_j$ (resp.\ $2_j$). 
See Figure~\ref{fig_cdp} again. 
We denote by $\overline{0_j 3_j}$ the edge in the piece $j$ 
joining vertices $0_j$ and $3_j$, and by $e_0$ the arc in $M_{n,k}$ 
(or the quotient space of $\mathcal{CP}_n$) obtained from 
gluing all $\overline{0_j 3_j}$'s. 

Let $D_{2 \, n}$ be the dihedral group of order $2 \, n$ generated by $r$ and $t$
(see Figure~\ref{fig_sym} again);
$D_{2 \, n} = \left< r, t \left| \, r^n = t^2 = (tr)^2 = 1 \right. \right>$.
Then we can easily see that, for all $n$ and $k$, 
$D_{2 \, n}$ is a subgroup of $\Aut \left( \mathcal{D}_{n,k} \right)$. 
So, without loss of generality, we only consider isomorphisms 
mapping $0_0$ to one of $0_0$, $1_0$, $0_1$, $1_1$. 
Since pieces of $\mathcal{CD}_{n,k}$ are not regular
tetrahedra, there are just the following  eight possibilities for such isomorphisms: 
\begin{align*}
\varphi_0 \colon (0_0, 1_0, 2_0, 3_0) &
\longmapsto (0_0, 1_0, 2_0, 3_0) \, , &
\varphi_1 \colon (0_0, 1_0, 2_0, 3_0) &
\longmapsto (0_1, 1_1, 2_1, 3_1) \, , \\
\varphi_2 \colon (0_0, 1_0, 2_0, 3_0) &
\longmapsto (1_0, 0_0, 3_0, 2_0) \, , &
\varphi_3 \colon (0_0, 1_0, 2_0, 3_0) &
\longmapsto (1_1, 0_1, 3_1, 2_1) \, , \\
\varphi_4 \colon (0_0, 1_0, 2_0, 3_0) &
\longmapsto (0_0, 3_0, 2_0, 1_0) \, , &
\varphi_5 \colon (0_0, 1_0, 2_0, 3_0) &
\longmapsto (0_1, 3_1, 2_1, 1_1) \, , \\
\varphi_6 \colon (0_0, 1_0, 2_0, 3_0) &
\longmapsto (1_0, 2_0, 3_0, 0_0) \, , &
\varphi_7 \colon (0_0, 1_0, 2_0, 3_0) &
\longmapsto (1_1, 2_1, 3_1, 0_1) \, .
\end{align*}
The isomorphism $\varphi_0$ induces the identity of 
$\Aut \left( \mathcal{D}_{n,k} \right)$, 
and it is easy to see that 
$\varphi_1 \left( \mathcal{D}_{n,k} \right) = \mathcal{D}_{n,n-k-1}$. 
So $M_{n,k} \cong M_{n',k'}$ if $n'=n$ and $k'=k$ or $n-k-1$. 
Furthermore $\varphi_1 \in \Aut \left( \mathcal{D}_{n,k} \right)$ 
if and only if $n-k-1=k$, i.e., $n$ is odd and $k = (n-1)/2$. 
%
%
%
%
%
%
%
%
\subsection*{Case 1 ($n \not\equiv 0 \bmod 3$).}
In this case, as we saw in the introduction,
all internal edges of $\mathcal{P}_n$ are glued together.
So the arc in $M_{n,k}$ obtained from the internal edges of 
$\mathcal{P}_n$ is surrounded by $6 \, n$ wedges. 
On the other hand, the arc $e_0$ in $M_{n,k}$ is surrounded by $2 \, n$ wedges. 
Thus any isomorphism must maps $\overline{0_i 3_i}$ to $\overline{0_j 3_j}$,
which means that
$\varphi_2$, $\varphi_3$, $\varphi_4$, $\varphi_5$, $\varphi_6$ and $\varphi_7$
are not combinatorial isomorphisms of $\mathcal{CD}_{n,k}$.

We first suppose that $n$ is odd and $k=(n-1)/2$.
Then $\varphi_0, \varphi_1 \in \Aut \left( \mathcal{D}_{n,k} \right)$
and we have the following short exact sequence: 
\begin{equation*}
1 \longrightarrow D_{2 \, n} \stackrel{\iota}{\longrightarrow}
\Aut \left( \mathcal{D}_{n,k} \right) \stackrel{\psi}{\longrightarrow}
{\Bbb Z}_2 \longrightarrow 1 \, ,
\end{equation*}
where $\iota$ is the inclusion, 
and $\psi$ is the homomorphism defined as follows: 
for any $\varphi \in \Aut \left( \mathcal{D}_{n,k} \right)$, 
\begin{equation*}
\psi (\varphi) := 
\begin{cases}
0& 
\text{if $\overline{0_0 3_0}$ maps to $\overline{0_i 3_i}$ for some even number $i$,}\\
1& 
\text{if $\overline{0_0 3_0}$ maps to $\overline{0_i 3_i}$ for some odd number $i$.}
\end{cases}
\end{equation*}
Let $u$ be $\varphi_1$, then 
${\Bbb Z}_2$ is presented as 
$\left< u \left| \, u^2 = 1 \right. \right>$ 
and is a subgroup of $\Aut \left( \mathcal{D}_{n,k} \right)$. 
So the sequence above splits. Thus we obtain 
$\Aut \left( \mathcal{D}_{n,k} \right)$ as follows: 
\begin{equation*}
\Aut \left( \mathcal{D}_{n,k} \right) \cong 
\left< t, u \left| \, t^2 = u^2 = (ut)^{2 \, n} = 1 \right. \right> \, .
\end{equation*}
We here note that $r=utut$. 

We secondly consider the other case, i.e., $n$ is even or $k \not= (n-1)/2$. 
In this case $u=\varphi_1 \not\in \Aut \left( \mathcal{D}_{n,k} \right)$. 
So $\psi$ becomes the trivial homomorphism and 
$\Aut \left( \mathcal{D}_{n,k} \right) \cong D_{2 \, n}$. 
%
%
%
%
%
%
%
%
\subsection*{Case 2 ($n \equiv 0 \bmod 3$).}
In this case, as mentioned in the introduction, 
the internal edges of $\mathcal{P}_n$ are
identified in three groups of $4 \, n / 3$ each. 
We denote by $e_1$ (resp.\ $e_2$, $e_3$) the arc in $M_{n,k}$ 
(or the quotient space of $\mathcal{CP}_n$) 
including $a_0$ (resp.\ $a_1$, $a_2$). 
This case is separated into two subcases.
%
%
%
%
%
%
\subsubsection*{Subcase~2.1 ($n \equiv 0$ {\rm and} $k \not\equiv 1 \bmod 3$).}
Though the arcs $e_0$, $e_1$, $e_2$ and $e_3$ are surrounded by $2 \, n$ wedges,
$e_0$ is done by $2 \, n$ different pieces, 
while each of $e_1$, $e_2$ and $e_3$ is done by $4 \, n/3$ different ones. 
So, in the same fashion as Case 1, 
$\varphi_2$, $\varphi_3$, $\varphi_4$, 
$\varphi_5$, $\varphi_6$ and $\varphi_7$ 
are not combinatorial isomorphisms of $\mathcal{CD}_{n,k}$.
Furthermore, there is no pair $(n,k)$ satisfying 
$n \equiv 0$ and $k \not\equiv 1$ together with $n$ is odd and $k = (n-1)/2$. 
This means that $\varphi_1 \not\in \Aut \left( \mathcal{D}_{n,k} \right)$, 
and thus we have $\Aut \left( \mathcal{D}_{n,k} \right) \cong D_{2 \, n}$.
%
%
%
%
%
%
\subsubsection*{Subcase~2.2 ($n \equiv 0$ {\rm and} $k \equiv 1 \bmod 3$).}
In this case $n = 3 \, m$ for some $m \geq 2$ and 
$k = 3 \, l + 1$ for some $0 \leq l \leq m-1$. 
Then all the arcs $e_0$, $e_1$, $e_2$ and $e_3$ are 
surrounded by $2 \, n$ different pieces. 

\begin{figure}[ht]
	\begin{center}
	\includegraphics[width=100mm, clip]{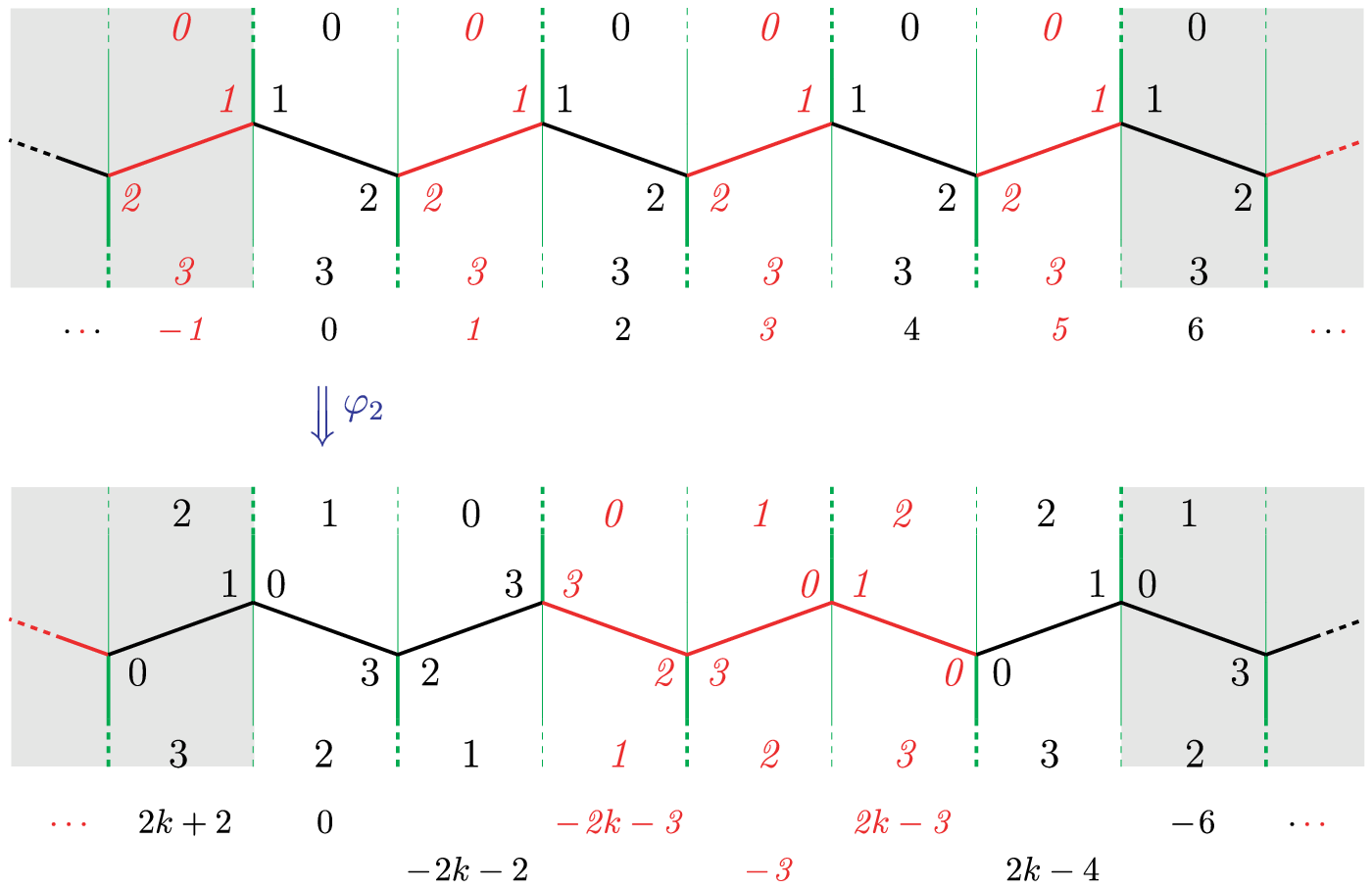}
	\end{center}
	\caption{A part of the side faces of 
		$\mathcal{CD}_{n,k}$ and 
		$\varphi_2 \left( \mathcal{CD}_{n,k} \right)$}
	\label{fig_phi21}
\end{figure}

We first show that 
$\varphi_2 \in \Aut \left( \mathcal{D}_{n,k} \right)$. 
A part of the side face of 
$\varphi_2 \left( \mathcal{CD}_{n,k} \right)$ 
looks as shown in Figure~\ref{fig_phi21}.
Since $k \equiv 1$,
the identification rule of step $k$ in terms of $\mathcal{CD}_{n,k}$
induces the one in terms of $\varphi_2 \left( \mathcal{CD}_{n,k} \right)$.
Thus we have $\varphi_2 \in \Aut \left( \mathcal{D}_{n,k} \right)$.
%
%
%
%
%
%

Now we have the following equalities:
$\varphi_3 = \varphi_1 \circ \varphi_2$,
$\varphi_4 = \varphi_1 \circ \varphi_5$,
$\varphi_5 = \varphi_2 \circ r^{-k-1}$,
$\varphi_6 = \varphi_2 \circ \varphi_4$,
$\varphi_7 = \varphi_5 \circ t$.
So all the possibilities for $\varphi_i \left( \mathcal{D}_{n,k} \right)$
are nothing but $\mathcal{D}_{n,k}$ or $\mathcal{D}_{n,n-k-1}$.
Thus we have proved \eqref{it_class}, i.e., 
$M_{n,k}$ and $M_{n',k'}$ are isometric if and only if
$n' = n$ and\/ $k' = k$ or\/ $n - k - 1$.

To finish the proof in this case,
we first consider the case where $n$ is even or $k \not= (n-1)/2$. 
In this case, as we saw before, 
$\varphi_1 \not\in \Aut \left( \mathcal{D}_{n,k} \right)$. 
So $\Aut \left( \mathcal{D}_{n,k} \right)$ is generated by 
$r$, $t$ and $s := \varphi_2$. 
In $M_{n,k}$ each element of $\Aut \left( \mathcal{D}_{n,k} \right)$
gives rise to a permutation of arcs $e_0$, $e_1$, $e_2$ and $e_3$. 
So there is a homomorphism, say $\psi$, from $\Aut \left( \mathcal{D}_{n,k} \right)$ 
to ${\frak S}_4$, the symmetric group of order $4$ acting on indices of arcs. 
Then we have $\psi(r) = (1 \, 2 \, 3)$, $\psi(s) = (0 \, 2)$ and $\psi(t) = (1 \, 3)$, 
and thus three elements $\psi(r)$, $\psi(s)$ and $\psi(t)$ generate ${\frak S}_4$ as follows: 
\begin{equation*}
{\frak S}_4 = \left< \psi(r), \psi(s), \psi(t) \left| \, 
\begin{array}{l}
\psi(r)^3 = \psi(s)^2 = \psi(t)^2 = (\psi(s) \psi(t) \psi(r))^3\\
= (\psi(s) \psi(t))^2 = (\psi(t) \psi(r))^2 = 1 
\end{array} \right. \right> \, .
\end{equation*}
So we have the following short exact sequence: 
\begin{equation*}
1 \longrightarrow \ker \psi \stackrel{\iota}{\longrightarrow}
\Aut \left( \mathcal{D}_{n,k} \right) \stackrel{\psi}{\longrightarrow}
{\frak S}_4 \longrightarrow 1 \, , 
\end{equation*}
where $\iota$ is the inclusion and $\ker \psi$ 
is isomorphic to ${\Bbb Z}_m$, generated by $r^3$. 
We here note that this sequence does not split.
Using combinatorial group theory  (see, for example, \cite[\S~10.2]{jo}), 
we obtain $\Aut \left( \mathcal{D}_{n,k} \right)$ as follows: 
\begin{equation*}
\Aut \left( \mathcal{D}_{n,k} \right) \cong 
\left< r, s, t \left| \, 
\begin{array}{l}
r^{3 \, m} = s^2 = t^2 = (tr)^2 = (st)^2 = 1,\\
sr^3 = r^3s, (str)^3 = r^{3 \left( m - 2 \, l - 2 \right)} 
\end{array} \right. \right> \left( =: G \right) \, .
\end{equation*}

We finally consider the case where $n$ is odd and $k = (n-1)/2$, i.e., 
$m$ is odd and $l = (m-1)/2$. 
In this case $\Aut \left( \mathcal{D}_{n,k} \right)$ 
is generated by $G$ and ${\Bbb Z}_2 \cong \left< u \left| \, u^2 = 1 \right. \right>$
because $u=\varphi_1 \in \Aut \left( \mathcal{D}_{n,k} \right)$.
Then $G$ is a normal subgroup of 
$\Aut \left( \mathcal{D}_{n,k} \right)$ since $usu = tsrs \in G$.
Furthermore $\Aut \left( \mathcal{D}_{n,k} \right)$
can be regarded as a semi-direct product of $G$ and ${\Bbb Z}_2$
since $G \cap {\Bbb Z}_2 = \left\{ 1 \right\}$.
Thus we obtain a presentation of $\Aut \left( \mathcal{D}_{n,k} \right)$ as follows:
\begin{equation*}
\Aut \left( \mathcal{D}_{n,k} \right) \cong 
\left< s, t, u \left| \,  
s^2 = t^2 = u^2 = (st)^2 = (ut)^6 = 1, susus = tutut \right. \right> \, .
\end{equation*}

We have thus finished the proof of the theorem. \qed


\providecommand{\bysame}{\leavevmode\hbox to3em{\hrulefill}\thinspace}

\end{document}